\documentclass[12pt]{article}%
\usepackage{amsfonts}
\usepackage{sw20bams}%
\usepackage{amsmath}%
\usepackage{amssymb}%
\usepackage{graphicx}
\providecommand{\U}[1]{\protect\rule{.1in}{.1in}}
\begin{document}

\title{Expected $n$-Step Product for Gaussian Tours}
\author{Steven Finch}
\date{December 17, 2015}
\maketitle

\begin{abstract}
Supplements to Mehta \&\ Normand (1997) are given, with regard to integrals
involving Euclidean distances between $n+1$ random points in $d$-dimensional
space, each visited once.

\end{abstract}

\footnotetext{Copyright \copyright \ 2015 by Steven R. Finch. All rights
reserved.}Let $\vec{r}_{1}$, $\vec{r}_{2}$, \ldots, $\vec{r}_{n+1}$ be
independent random Gaussian points in $\mathbb{R}^{d}$, all of which have mean
vector zero and covariance matrix identity. \ What can be said about%
\[
\mu_{d,n}=\operatorname{E}\left(  \left\vert \vec{r}_{2}-\vec{r}%
_{1}\right\vert \cdot\left\vert \vec{r}_{3}-\vec{r}_{2}\right\vert
\cdot\left\vert \vec{r}_{4}-\vec{r}_{3}\right\vert \cdots\left\vert \vec
{r}_{n}-\vec{r}_{n-1}\right\vert \cdot\left\vert \vec{r}_{n+1}-\vec{r}%
_{n}\right\vert \right)  ,
\]%
\[
\nu_{d,n}=\operatorname{E}\left(  \left\vert \vec{r}_{2}-\vec{r}%
_{1}\right\vert \cdot\left\vert \vec{r}_{3}-\vec{r}_{2}\right\vert
\cdot\left\vert \vec{r}_{4}-\vec{r}_{3}\right\vert \cdots\left\vert \vec
{r}_{n}-\vec{r}_{n-1}\right\vert \cdot\left\vert \vec{r}_{1}-\vec{r}%
_{n}\right\vert \right)
\]
? \ In words, what is the expected product of lengths of an $n$-step
open/closed tour of the points? \ \textquotedblleft Open\textquotedblright%
\ means that each point is uniquely visited; \textquotedblleft
closed\textquotedblright\ means the same except $\vec{r}_{n+1}$ is replaced by
$\vec{r}_{1}$. \ We do not assume anything about the ordering of the points,
hence if $n>d+1$ the closed tour need not lie on a convex hull boundary.

Mehta \&\ Normand \cite{MN-StpPrdTrs} computed%
\[%
\begin{array}
[c]{ccc}%
\mu_{3,1}=\dfrac{4}{\sqrt{\pi}}=2.256758..., &  & \mu_{3,2}=2+\dfrac{6\sqrt
{3}}{\pi}=5.307973...,
\end{array}
\]%
\[
\mu_{3,3}=\frac{238}{3\sqrt{\pi}}+\frac{56\sqrt{2}}{3\pi^{3/2}}-\frac{216}%
{\pi^{3/2}}\arctan\left(  \sqrt{2}\right)  =12.442385...,
\]%
\begin{align*}
\mu_{3,4}  & =\frac{232}{45}-\frac{3140}{9\pi}+\frac{56}{\sqrt{3}\pi}%
+\frac{260\sqrt{5}}{9\pi^{2}}+\frac{912}{\pi^{2}}\arctan\left(  \sqrt
{5}\right)  +\frac{224}{\sqrt{3}\pi^{2}}\arctan\left(  \sqrt{\frac{5}{3}%
}\right) \\
& =29.174181...,
\end{align*}%
\[
\nu_{3,2}=6
\]
using a multipole expansion (via spherical harmonics) and solution of Fredholm
integral equations (via eigenfunction analysis). \ No exact pattern for either
$\mu_{3,n}$ or $\nu_{3,n}$ is observed. \ It is tempting to believe that the
difficulties encountered when $d=3$ might be somehow circumvented when $d=1$
or $d=2$. \ This turns out to be false. \ On a line,
\[%
\begin{array}
[c]{ccc}%
\mu_{1,1}=\dfrac{2}{\sqrt{\pi}}=1.128379..., &  & \mu_{1,2}=\dfrac{1}%
{3}+\dfrac{2\sqrt{3}}{\pi}=1.435991...,
\end{array}
\]%
\[
\mu_{1,3}=\frac{5}{\sqrt{\pi}}+\frac{4\sqrt{2}}{\pi^{3/2}}-\frac{12}{\pi
^{3/2}}\arctan\left(  \sqrt{2}\right)  =1.778095...,
\]%
\[
\mu_{1,4}=\frac{2}{15}+\frac{4\sqrt{5}}{\pi^{2}}+\dfrac{8}{\pi^{2}}%
\arctan\left(  \frac{\sqrt{5}}{7}\right)  +\dfrac{8\sqrt{3}}{\pi^{2}}%
\arctan\left(  \sqrt{\frac{3}{5}}\right)  =2.215483...,
\]%
\[%
\begin{array}
[c]{ccccc}%
\nu_{1,2}=2, &  & \nu_{1,3}=\dfrac{3}{\sqrt{\pi}}=1.692568..., &  & \nu
_{1,4}=\dfrac{2}{3}-\dfrac{8}{\pi}+\dfrac{8\sqrt{3}}{\pi}=2.530818....
\end{array}
\]
The formula for $\mu_{1,3}$ in \cite{MN-StpPrdTrs} contains typographical
errors and that for $\mu_{1,4}$ does not appear at all; we give proofs in
Section 2. \ In the plane,
\[
\mu_{2,1}=\sqrt{\pi}=1.772453...,
\]%
\[
\mu_{2,2}=4E\left(  \dfrac{1}{2}\right)  -\dfrac{3}{2}K\left(  \dfrac{1}%
{2}\right)  =3.341223..,
\]%
\[
\nu_{2,2}=4
\]
where
\[%
\begin{array}
[c]{c}%
K(\xi)=%
{\displaystyle\int\limits_{0}^{\pi/2}}
\dfrac{1}{\sqrt{1-\xi^{2}\sin(\theta)^{2}}}\,d\theta=%
{\displaystyle\int\limits_{0}^{1}}
\dfrac{1}{\sqrt{(1-t^{2})(1-\xi^{2}t^{2})}}\,dt,\\
E(\xi)=%
{\displaystyle\int\limits_{0}^{\pi/2}}
\sqrt{1-\xi^{2}\sin(\theta)^{2}}\,d\theta=%
{\displaystyle\int\limits_{0}^{1}}
\sqrt{\dfrac{1-\xi^{2}t^{2}}{1-t^{2}}}\,dt
\end{array}
\]
are complete elliptic integrals of the first and second kind
\cite{M1-StpPrdTrs, M2-StpPrdTrs}. \ Also \cite{CG-StpPrdTrs, Fi-StpPrdTrs},%
\begin{align*}
\nu_{2,3}  & =\dfrac{4}{3\pi}%
{\displaystyle\int\limits_{0}^{\infty}}
{\displaystyle\int\limits_{0}^{x}}
{\displaystyle\int\limits_{x-y}^{x+y}}
\dfrac{x^{2}y^{2}z^{2}}{\sqrt{(x+y+z)(-x+y+z)(x-y+z)(x+y-z)}}\cdot\\
& \exp\left(  -\dfrac{1}{6}\left(  x^{2}+y^{2}+z^{2}\right)  \right)
dz\,dy\,dx\\
& \approx6.359,
\end{align*}%
\[
\mu_{2,3}=\lim_{\rho\rightarrow-\frac{1}{2}^{+}}F(\rho)\approx6.25
\]
and we examine the latter in Section 3. \ An expression for $F(\rho)$ and
explanation of its significance are forthcoming. \ Finally, returning to
three-space \cite{CG-StpPrdTrs, Fi-StpPrdTrs},%
\[
\nu_{3,3}=\dfrac{2\sqrt{3}}{9\pi}%
{\displaystyle\int\limits_{0}^{\infty}}
{\displaystyle\int\limits_{0}^{x}}
{\displaystyle\int\limits_{x-y}^{x+y}}
\,x^{2}y^{2}z^{2}\exp\left(  -\dfrac{1}{6}\left(  x^{2}+y^{2}+z^{2}\right)
\right)  dz\,dy\,dx\approx12.708.
\]
Accurate numerical values for $\mu_{2,4}$, $\nu_{2,4}$, $\nu_{3,4}$ (obtained
by procedures other than Monte Carlo simulation) remain at large.

\section{One-Space}

Define correlation coefficients%
\begin{align*}
\rho_{ij}  & =\frac{1}{2}\operatorname*{E}\left(  (r_{i+1}-r_{i}%
)(r_{j+1}-r_{j})\right) \\
& =\left\{
\begin{array}
[c]{ccc}%
1 &  & \text{if }\left\vert i-j\right\vert =0,\\
-1/2 &  & \text{if }\left\vert i-j\right\vert =1,\\
0 &  & \text{otherwise}%
\end{array}
\right.
\end{align*}
and partial correlation coefficients
\[
\rho_{ij\cdot k}=\frac{\rho_{ij}-\rho_{ik}\rho_{jk}}{\sqrt{\left(  1-\rho
_{ik}^{2}\right)  \left(  1-\rho_{jk}^{2}\right)  }},
\]%
\[
\rho_{ij\cdot k\ell}=\frac{\rho_{ij\cdot k}-\rho_{i\ell\cdot k}\rho
_{j\ell\cdot k}}{\sqrt{\left(  1-\rho_{i\ell\cdot k}^{2}\right)  \left(
1-\rho_{j\ell\cdot k}^{2}\right)  }}.
\]
(In statistics \cite{BS-StpPrdTrs}, $\rho_{ij\cdot k\ell}=0$ would suggest
that $r_{i+1}-r_{i}$ and $r_{j+1}-r_{j}$ are conditionally independent, given
$r_{k+1}-r_{k}$ and $r_{\ell+1}-r_{\ell}$. \ More precisely, $\rho_{ij\cdot
k\ell}$ is the correlation between the residuals for $i$ and $j$ resulting
from a linear regression of $i$ with $k,\ell$ and of $j$ with $k,\ell$.) \ Let
$R_{n}$ be the determinant of $\left(  \rho_{ij}\right)  _{1\leq i,j\leq n}$.
\ It follows from general formulas in \cite{N1-StpPrdTrs, N2-StpPrdTrs,
N3-StpPrdTrs} that%
\[
\mu_{1,2}=\frac{4}{\pi}\left(  \sqrt{R_{2}}+\rho_{12}\arcsin(\rho
_{12})\right)  ,
\]%
\begin{align*}
\mu_{1,3}  & =\frac{8}{\pi^{3/2}}\left[  \sqrt{R_{3}}+\left(  \rho_{12}%
+\rho_{13}\rho_{23}\right)  \arcsin(\rho_{12\cdot3})+\right. \\
& \left.  \left(  \rho_{13}+\rho_{12}\rho_{23}\right)  \arcsin(\rho_{13\cdot
2})+\left(  \rho_{23}+\rho_{12}\rho_{13}\right)  \arcsin(\rho_{23\cdot
1})\right]  ,
\end{align*}%
\begin{align*}
\mu_{1,4}  & =\frac{16}{\pi^{2}}\left[  \sqrt{R_{4}}+\sqrt{1-\rho_{12}^{2}%
}\left(  \rho_{34}+\rho_{13}\rho_{14}+\rho_{23}\rho_{24}\right)  \arcsin
(\rho_{34\cdot12})+\right. \\
& \sqrt{1-\rho_{13}^{2}}\left(  \rho_{24}+\rho_{12}\rho_{14}+\rho_{23}%
\rho_{34}\right)  \arcsin(\rho_{24\cdot13})+\\
& \sqrt{1-\rho_{14}^{2}}\left(  \rho_{23}+\rho_{12}\rho_{13}+\rho_{24}%
\rho_{34}\right)  \arcsin(\rho_{23\cdot14})+\\
& \sqrt{1-\rho_{23}^{2}}\left(  \rho_{14}+\rho_{12}\rho_{24}+\rho_{13}%
\rho_{34}\right)  \arcsin(\rho_{14\cdot23})+\\
& \sqrt{1-\rho_{24}^{2}}\left(  \rho_{13}+\rho_{12}\rho_{23}+\rho_{14}%
\rho_{34}\right)  \arcsin(\rho_{13\cdot24})+\\
& \left.  \sqrt{1-\rho_{34}^{2}}\left(  \rho_{12}+\rho_{13}\rho_{23}+\rho
_{14}\rho_{24}\right)  \arcsin(\rho_{12\cdot34})\right]  +\\
& 4\left(  \rho_{12}\rho_{34}+\rho_{13}\rho_{24}+\rho_{14}\rho_{23}\right)
\gamma
\end{align*}
where
\[
\gamma=\operatorname{E}\left[  \operatorname*{sgn}\left(  (r_{2}-r_{1}%
)(r_{3}-r_{2})(r_{4}-r_{3})(r_{5}-r_{4})\right)  \right]
\]
and $\operatorname*{sgn}(x)=1$ when $x\geq0$; $\operatorname*{sgn}(x)=-1$ when
$x<0$. \ The additional factor $\gamma$ could be troublesome were it not known
that the orthant probability
\[
\operatorname*{P}\left\{  r_{2}-r_{1}>0,r_{3}-r_{2}>0,r_{4}-r_{3}%
>0,r_{5}-r_{4}>0\right\}  =\frac{1}{16}+\frac{1}{8\pi}%
{\displaystyle\sum\limits_{i<j}}
\arcsin(\rho_{ij})+\frac{\gamma}{16}
\]
is equal to $1/120$ \cite{P-StpPrdTrs, M-StpPrdTrs, C-StpPrdTrs,
SA-StpPrdTrs}. \ From this, we deduce that $\gamma=2/15$ and thus the value
for $\mu_{1,4}$ is confirmed.

\section{Two-Space}

Thinking of vectors as complex numbers:%
\[%
\begin{array}
[c]{ccc}%
\vec{r}_{2}-\vec{r}_{1}=z_{1}=\left(  x_{2}-x_{1}\right)  +i\left(
y_{2}-y_{1}\right)  , &  & \vec{r}_{3}-\vec{r}_{2}=z_{2}=\left(  x_{3}%
-x_{2}\right)  +i\left(  y_{3}-y_{2}\right)
\end{array}
\]
we have%
\[%
\begin{array}
[c]{ccc}%
\bar{z}_{1}=\left(  x_{2}-x_{1}\right)  -i\left(  y_{2}-y_{1}\right)  , &  &
\bar{z}_{2}=\left(  x_{3}-x_{2}\right)  -i\left(  y_{3}-y_{2}\right)
\end{array}
\]
and thus%
\[
\psi_{11}=\operatorname*{E}\left(  z_{1}\bar{z}_{1}\right)  =\operatorname*{E}%
\left(  x_{2}^{2}+x_{1}^{2}+y_{2}^{2}+y_{1}^{2}\right)  =4,
\]%
\[
\psi_{12}=\operatorname*{E}\left(  z_{1}\bar{z}_{2}\right)  =\operatorname*{E}%
\left(  -x_{2}^{2}-y_{2}^{2}\right)  =-2.
\]
Writing%
\[
\Psi=\left(
\begin{array}
[c]{cc}%
\psi_{11} & \psi_{12}\\
\bar{\psi}_{12} & \psi_{22}%
\end{array}
\right)  =\left(
\begin{array}
[c]{cc}%
4 & -2\\
-2 & 4
\end{array}
\right)  ,
\]
we obtain inverse covariance matrix%
\[
\Phi=\left(
\begin{array}
[c]{cc}%
\varphi_{11} & \varphi_{12}\\
\bar{\varphi}_{12} & \varphi_{22}%
\end{array}
\right)  =\Psi^{-1}=\left(
\begin{array}
[c]{cc}%
1/3 & 1/6\\
1/6 & 1/3
\end{array}
\right)
\]
with determinant $\Delta=1/12$. \ The joint density of $a=\left\vert
z_{1}\right\vert $, $b=\left\vert z_{2}\right\vert $ is \cite{M1-StpPrdTrs,
M2-StpPrdTrs}%
\[
4\Delta\,a\,b\exp\left(  -\varphi_{11}a^{2}-\varphi_{22}b^{2}\right)
I_{0}\left(  2a\,b\left\vert \varphi_{12}\right\vert \right)
\]
where $I_{k}$ is the $k^{\text{th}}$ modified Bessel function of the first
kind, and therefore%
\[
\mu_{2,2}=\frac{1}{3}%
{\displaystyle\int\limits_{0}^{\infty}}
{\displaystyle\int\limits_{0}^{\infty}}
a^{2}b^{2}\exp\left(  -\dfrac{1}{3}\left(  a^{2}+b^{2}\right)  \right)
I_{0}\left(  \frac{a\,b}{3}\right)  db\,da=4E\left(  \dfrac{1}{2}\right)
-\dfrac{3}{2}K\left(  \dfrac{1}{2}\right)  .
\]

The joint density of $a=\left\vert z_{1}\right\vert $, $b=\left\vert
z_{2}\right\vert $, $c=\left\vert z_{3}\right\vert $ is more complicated. \ It
is useful to introduce a parameter $\rho$ for which%
\[
\Psi(\rho)=4\left(
\begin{array}
[c]{ccc}%
1 & \rho & 0\\
\rho & 1 & \rho\\
0 & \rho & 1
\end{array}
\right)  \rightarrow\left(
\begin{array}
[c]{ccc}%
4 & -2 & 0\\
-2 & 4 & -2\\
0 & -2 & 4
\end{array}
\right)
\]
as $\rho\rightarrow-1/2$. \ The inverse covariance matrix is%
\[
\Phi(\rho)=\frac{1}{4\left(  1-2\rho^{2}\right)  }\left(
\begin{array}
[c]{ccc}%
1-\rho^{2} & -\rho & \rho^{2}\\
-\rho & 1 & -\rho\\
\rho^{2} & -\rho & 1-\rho^{2}%
\end{array}
\right)
\]
and possesses determinant%
\[
\Delta(\rho)=\frac{1}{64\left(  1-2\rho^{2}\right)  }.
\]
Restrict $-1/\sqrt{2}<\rho<0$ so that all $\varphi_{ij}$ and $\Delta$ are
positive. Let $\varepsilon_{0}=1$, $\varepsilon_{k}=2$ for $k\geq1$. \ The
joint density $f(a,b,c)$ is given by an infinite series \cite{M1-StpPrdTrs,
M2-StpPrdTrs}%
\begin{align*}
& 8\Delta\,a\,b\,c\exp\left(  -\varphi_{11}a^{2}-\varphi_{22}b^{2}%
-\varphi_{33}c^{2}\right)
{\displaystyle\sum\limits_{k=0}^{\infty}}
\varepsilon_{k}(-1)^{k}\cdot\\
& I_{k}\left(  2a\,b\left\vert \varphi_{12}\right\vert \right)  I_{k}\left(
2b\,c\left\vert \varphi_{23}\right\vert \right)  I_{k}\left(  2a\,c\left\vert
\varphi_{13}\right\vert \right) \\
& \smallskip\\
& =\frac{a\,b\,c}{8\left(  1-2\rho^{2}\right)  }\exp\left(  -\frac{\left(
1-\rho^{2}\right)  a^{2}+b^{2}+\left(  1-\rho^{2}\right)  c^{2}}{4\left(
1-2\rho^{2}\right)  }\right)
{\displaystyle\sum\limits_{k=0}^{\infty}}
\varepsilon_{k}(-1)^{k}\cdot\\
& I_{k}\left(  \frac{-\rho\,a\,b}{2\left(  1-2\rho^{2}\right)  }\right)
I_{k}\left(  \frac{-\rho\,b\,c}{2\left(  1-2\rho^{2}\right)  }\right)
I_{k}\left(  \frac{\rho^{2}a\,c}{2\left(  1-2\rho^{2}\right)  }\right)
\end{align*}
which is evidently divergent when $\rho=-1/2$. \ Convergence seems to
stabilize for $\rho$ slightly to the right of $-1/2$. We therefore examine%
\[
\mu_{2,3}=\lim_{\rho\rightarrow-\frac{1}{2}^{+}}%
{\displaystyle\int\limits_{0}^{\infty}}
{\displaystyle\int\limits_{0}^{\infty}}
{\displaystyle\int\limits_{0}^{\infty}}
a\,b\,c\,f(a,b,c)dc\,db\,da\approx6.25,
\]
yielding a result consistent with (but unfortunately not improving upon)
computer simulation.

A\ higher-dimensional analog of $f(a,b,c)$ is exhibited in \cite{BM-StpPrdTrs,
CT-StpPrdTrs} but requires that the $4\times4$ inverse covariance matrix have
corner entry $\varphi_{14}=0$. \ This is not true in our scenario. \ More
relevant discussion is found in \cite{NK1-StpPrdTrs, NK2-StpPrdTrs,
LW1-StpPrdTrs, LW2-StpPrdTrs}.


\begin{thebibliography}{99}                                                                                               %
\bibitem {MN-StpPrdTrs}M. L. Mehta and J.-M. Normand, On some definite
multiple integrals, \textit{J. Phys. A} 30 (1997) 8671--8684; MR1619533 (99g:33059).

\bibitem {M1-StpPrdTrs}K. S. Miller, Complex Gaussian processes, \textit{SIAM
Rev.} 11 (1969) 544--567; MR0258109 (41 \#2756).

\bibitem {M2-StpPrdTrs}K. S. Miller, \textit{Complex Stochastic Processes. An
Introduction to Theory and Application}, Addison-Wesley, 1974, pp. 86--100;
MR0368118 (51 \#4360).

\bibitem {CG-StpPrdTrs}P. Clifford and N. J. B. Green, Distances in Gaussian
point sets, \textit{Math. Proc. Cambridge Philos. Soc.} 97 (1985) 515--524;
MR0778687 (86i:62091).

\bibitem {Fi-StpPrdTrs}S. R. Finch, Random triangles, unpublished note (2010),
http://www.people.fas.harvard.edu/\symbol{126}sfinch/.

\bibitem {BS-StpPrdTrs}K. Baba, R. Shibata and M. Sibuya, Partial correlation
and conditional correlation as measures of conditional independence,
\textit{Austral. New Zealand J. Statist.} 46 (2004) 657--664; MR2115961 (2005k:62153).

\bibitem {N1-StpPrdTrs}S. Nabeya, Absolute moments in $2$-dimensional normal
distribution, \textit{Annals Inst. Statist. Math.} 3 (1951) 2--6; MR0045347 (13,570b).

\bibitem {N2-StpPrdTrs}S. Nabeya, Absolute moments in $3$-dimensional normal
distribution, \textit{Annals Inst. Statist. Math.} 4 (1952) 15--30; MR0052072 (14,569c).

\bibitem {N3-StpPrdTrs}S. Nabeya, Absolute and incomplete moments of the
multivariate normal distribution, \textit{Biometrika} 48 (1961) 77--84;
MR0126917 (23 \#A4211).

\bibitem {P-StpPrdTrs}R. L. Plackett, A reduction formula for normal
multivariate integrals, \textit{Biometrika} 41 (1954) 351--360; MR0065047 (16,377c).

\bibitem {M-StpPrdTrs}J. A. McFadden, Two expansions for the quadrivariate
normal integral, \textit{Biometrika} 47 (1960) 325--333; MR0119221 (22 \#9987).

\bibitem {C-StpPrdTrs}D. R. Childs, Reduction of the multivariate normal
integral to characteristic form, \textit{Biometrika} 54 (1967) 293--300;
MR0214177 (35 \#5028).

\bibitem {SA-StpPrdTrs}H.-J. Sun and C. Asano, On the normal orthant
probability with a tri-diagonal correlation matrix, \textit{Eng. Sci. Rep.}
\textit{Kyushu Univ. }12 (1990) 53--58.

\bibitem {BM-StpPrdTrs}L. E. Blumenson and K. S. Miller, Properties of
generalized Rayleigh distributions, \textit{Annals Math. Statist.} 34 (1963)
903--910; MR0150860 (27 \#846).

\bibitem {CT-StpPrdTrs}Y. Chen and C. Tellambura, Infinite series
representations of the trivariate and quadrivariate Rayleigh distribution and
their applications, \textit{IEEE Trans. Comm}. 53 (2005) 2092--2101.

\bibitem {NK1-StpPrdTrs}S. Nadarajah and S. Kotz, On the infinite series
representations for multivariate Rayleigh distributions, \textit{IEEE Trans.
Comm}. 55 (2007) 392--393.

\bibitem {NK2-StpPrdTrs}S. Nadarajah and S. Kotz, Comments on
\textquotedblleft A trivariate chi-squared distribution derived from the
complex Wishart distribution\textquotedblright, \textit{J. Multivariate Anal.}
99 (2008) 306--307; MR2432328.

\bibitem {LW1-StpPrdTrs}W. V. Li and A. Wei, Gaussian integrals involving
absolute value functions, \textit{High Dimensional Probability V}, Proc. 2008
Luminy conf., ed. C. Houdr\'{e}, V. Koltchinskii, D. M. Mason and M. Peligrad,
Inst. Math. Statist., 2009, pp. 43--59; MR2797939 (2012f:60059).

\bibitem {LW2-StpPrdTrs}W. V. Li and A. Wei, A Gaussian inequality for
expected absolute products, \textit{J. Theoret. Probab.} 25 (2012) 92--99; MR2886380.%

\begin{tabular}
[c]{lll}
& Steven Finch & \\
& Dept. of Statistics & \\
& Harvard University & \\
& Cambridge, MA, USA & \\
& \textit{steven\_finch@harvard.edu} &
\end{tabular}

\end{thebibliography}
\end{document}